\documentclass[11pt, amssymb]{amsart}
\usepackage{amssymb}
\usepackage[mathscr]{eucal}

\def\cD{{\mathcal D}}

\def\cR{{\mathscr R}}
\def\cG{{\mathcal G}}

\def\cT{{\mathcal T}}
\def\bR{{\mathbb R}}
\def\bZ{{\mathbb Z}}
\def\bQ{{\mathbb Q}}
\def\bC{{\mathbb C}}

\def\cT{{\mathscr T}}

\def\cC{{\mathscr C}}
\def\cD{{\mathscr D}}

\def\ni{\noindent}

\def\bP{{\bf P}}

\def\oz1{d{\overline z}^1}
\def\oz2{d{\overline z}^2}
\def\oz3{d{\overline z}^3}

\def\oI{\overline I}

\def\oz{\overline z}

\def\oIq1{\oI_1\cdots\oI_{q-1}}
\def\oIq2{\oI_1\cdots\oI_{q-2}}

\def\fv{\mathfrak{v}}
\def\fD{\mathfrak{D}}
\begin{document}
\vskip5mm

\centerline{\bf Addendum to} 

\vskip2mm 
   
\centerline{\bf Fake Projective Planes}
  \vskip1.mm
 
 \centerline{{\it Inventiones Math.}\,{\bf 168}, 321-370 (2007)}
\vskip1.5mm
 
 \centerline{By Gopal Prasad and Sai-Kee Yeung}

\vskip5mm
We will use the notions and notations introduced in [1]. 

\vskip 1mm

\ni{\bf A.1.}  As in [1], let $k$ be a totally real number field and $\ell$ be a totally complex quadratic extension of $k$. Let $h$ be a hermitian form on $\ell^3$, defined in terms of the nontrivial automorphism of $\ell/k$,  such that $h$ is indefinite exactly at one real place of $k$.   
In [1] we proved that if the fundamental group of a fake projective plane is an arithmetic subgroup of 
${\rm PU}(h)$, then the pair $(k,\ell)$ must be one of the following five: $\cC_1$, $\cC_8$, $\cC_{11}$, $\cC_{18}$, or $\cC_{21}$. Around the time [1] was submitted for publication, Tim Steger had shown that $(k,\ell)$ cannot be $\cC_8$, and together with 
Donald Cartwright he had shown that it cannot be $\cC_{21}$. Using long and very sophisticated computer-assisted group theoretic computations, Donald Cartwright and Tim Steger have now shown that $(k,\ell)$ cannot be any of the three remaining pairs $\cC_1$, $\cC_{11}$ and $\cC_{18}$.  This then proves that the fundamental group of a fake projective plane cannot be an arithmetic subgroup of ${\rm PU}(h)$.  Cartwright and Steger have also shown that there exists a rather unexpected smooth projective complex algebraic surface whose Euler-Poincar\'e characteristic is $3$ but the first Betti-number is nonzero (it is actually $2$). This surface  is uniformized by the complex 2-ball, and its fundamental group is a cocompact torsion-free arithmetic subgroup of ${\rm PU}(h)$, $h$ as above, with $(k,\ell) = \cC_{11}$ $=(\bQ(\sqrt{3}),\bQ(\zeta_{12}) $).  

\vskip1mm
In view of the above, we can assert now that there are exactly {\it twenty eight} classes of fake projective planes in a classification scheme which is finer than that used in [1], see A.2 and A.5 below. It has been shown  recently  by Cartwright and Steger using ingenious group theoretic computations that  the twenty eight classes of fake projective planes altogether contain {\it fifty} distinct
fake projective planes up to isometry with respect to the Poincar\'e metric.  Since each such fake projective plane as a Riemannian manifold supports
two distinct complex structures [KK, \S5], there are exactly {\it one hundred} fake  projective planes counted 
up to biholomorphism.
\vskip1mm

 Cartwright and Steger have  given explicit generators and relations for the fundamental group (which is a cocompact torsion-free arithmetic subgroup of ${\rm PU}(2,1)$) of each of the fake projective planes, determined their automorphism group, and computed their first homology with coefficients in $\bZ$. They have shown that the fundamental group of eight of the one hundred fake projective planes do not admit an embedding into ${\rm SU}(2,1)$ as a discrete subgroup (the canonical line bundle of these fake projective planes is not divisible by $3$ in their Picard group).  All these fake projective planes arise from the pair $\cC_{18}$. Cartwright and Steger have also shown that the quotient of six of the fake projective planes by a subgroup of order $3$ of the automorphism group is a simply connected singular surface.

\vskip1mm

\vskip2mm

\ni{\bf A.2.} In 1.3 of [1], for any nonarchimedean place $v$ of $k$ which ramifies in $\ell$, $P_v$ was chosen to be a maximal parahoric subgroup of $G(k_v)$ normalized by $\Pi$, and in 1.4, in the collection $(P_v)_{v\in V_f}$, the parahoric subgroup $P_v$ was assumed to be maximal for all  $v$ which ramify  in $\ell$. It appears to us, in retrospect,  that we should not have required $P_v$ to be maximal for $v$ ramifying in $\ell$. We therefore abandon this condition. We then get a finer classification of the fake projective planes. 
\vskip1mm

Let  $\cT$ be the set of nonarchimedean places $v$ of $k$ such that in the collection $(P_v)_{v\in V_f}$ under consideration, $P_v$ is not maximal, and also all those $v$ which are unramified in $\ell$ and $P_v$ is not hyperspecial. Let $\cT_0$ be the subset of $\cT$ consisting of places where the group $G$ is anisotropic. 
\vskip1mm

In Theorem 4.4 of [1], the assertion that $\cT =\cT_0$ needs to be modified for $\ell = \bQ(\sqrt{-1})$, $\bQ(\sqrt{-2})$, and $\bQ(\sqrt{-7})$. It can be easily seen, by simple computations, that for $\ell =\bQ(\sqrt{-1})$ and $\bQ(\sqrt{-2})$, the possibilities are $\cT= \cT_0$ and $\cT = \cT_0\cup\{2\}$. For $\ell = \bQ(\sqrt{-7})$ the possibilities for $\cT$ are $\cT =\cT_0 =\{2\}$, $\cT =\{2,3\}$, and $\cT=\{2,5\}$. 
\vskip1mm

In this paragraph $a$ is either $1$ or $2$. To  see the fake projective planes arising from $\ell = \bQ(\sqrt{-a})$, and  $\cT= \cT_0 \cup\{2\}$,  we recall from [1] that  $\cT_0 = \{p\}$, where $p$ is $5$ if $a =1$, and is $3$ for $a =2$. Let the $\bQ$-group $G$ be as in 5.7 of [1], and $\overline{G}$ be its adjoint group.   Note that $2$ is the only prime which ramifies in $\bQ(\sqrt{-a})$.  Now let $(P_q)$ be a coherent collection of maximal parahoric subgroups $P_q$ of $G(\bQ_q)$, we assume that for $q\ne 2,\, p$, $P_q$ is hyperspecial. We fix an Iwahori subgroup $I_2$ contained in $P_2$ and let $\Lambda = G(\bQ)\cap\prod_q P_q$, $\Lambda_I = \Lambda\cap I_2 $.  Then $\Lambda$, and so also $\Lambda_I$, is torsion-free (Lemma 5.6 of [1]). Since $[P_2:I_2] =3$, the strong approximation property implies that $[\Lambda:\Lambda_I] = 3$. As $\chi(\Lambda) = 3$, we obtain that $\chi(\Lambda_I ) = 9$.  Now let $\Gamma_I$ be the normalizer of $\Lambda_I$ in $G(\bR)$, and $\overline{\Lambda}_I$ and $\overline{\Gamma}_I$ be the images of $\Lambda_I$ and $\Gamma_I$ in $\overline{G}(\bR)$.  Then $[\Gamma_I:\Lambda_I] = 9$ (cf.\,5.4 of [1]), and hence, $[\overline{\Gamma}_I:\overline{\Lambda}_I] =3$. This implies that the orbifold Euler-Poincar\'e characteristic of $\overline{\Gamma}_I$ is $3$. Moreover, $H^1(\overline{\Lambda}_I,\bC)$, and hence also $H^1(\overline{\Gamma}_I,\bC)$ vanishes (Theorem 15.3.1 of [Ro]).  Thus if $\overline{\Gamma}_I$ is torsion-free, which indeed is the case, as can be seen by a suitable adaptation of the argument used to prove the first assertion of the proposition in A.6 below, $B/\overline{\Gamma}_I$ is a fake projective plane, and it is the unique plane belonging to the class associated to $\ell = \bQ(\sqrt{-a})$, for $a = 1,\,2$, and $\cT = \cT_0\cup \{2\}$.

 \vskip1mm
 
In view of the above, under the new, and finer, classification scheme, the fake projective planes arising from  $\ell =\bQ(\sqrt{-1})$ and $\bQ(\sqrt{-2})$ are placed in six rather than the original four classes.       
\vskip1mm

\vskip1mm

We note that the fact that $\overline{\Gamma}_I$ described above is torsion-free, and the proposition in A.6 below, have  also been proved by Cartwright and Steger using computer-assisted group theoretic computations. 
 
\vskip2mm

\ni{\bf A.3.} In this section, $\ell =\bQ(\sqrt{-7})$. As in 5.7 of [1], let $\cD$ be a cubic division algebra with center $\ell$ whose local invariants at the two places of $\ell$ lying over $p=2$ are nonzero and negative of each other, and whose local invariant at all the  other places of $\ell$ is zero. Let the $\bQ$-group $G$ be as in 5.5 of [1]  for an involution $\sigma$ of $\cD$ of  the second kind chosen so that $G(\bR)\cong{\rm{SU}}(2,1)$.  Then 
$\cT_0=\{2\}$. Let $C$ be the center of $G$ and $\overline{G}= G/C$ be the adjoint group of $G$. As $\ell$ does not contain a nontrivial cube-root of unity, $C(\bQ)$ is trivial.   
\vskip2mm

The following seven paragraphs  should be added at the end of 5.11 in [1]. 
\vskip2mm

Let $\cT =\{2,3\}$ or $\{2,5\}$. We fix a coherent collection $(P_q)$ of {\it maximal} parahoric subgroups $P_q$ of $G(\bQ_q)$. We assume that $P_q$ is hyperspecial if, and only if, $q\notin \cT\cup\{7\}$. Let $\Lambda =G(\bQ)\cap\prod_q P_q$, and let $\Gamma$ be its normalizer in $G(\bR)$.
Let $\overline\Gamma$ (resp.,\,$\overline\Lambda$) be the image of $\Gamma$ (resp.,\,$\Lambda$) in $\overline{G}(\bR)$. Then as $\#\cT_0 = 1$, $[\Gamma:\Lambda]=9$, see 5.4, and hence, $[\overline\Gamma :\overline\Lambda ]= 3$.  We note that $\overline\Gamma$ is contained in $\overline{G}(\bQ)$ (see, for example, [BP, Proposition 1.2]). Obviously, $\overline\Gamma$ normalizes $P_q$ for every $q$.

\vskip1mm

\ni (i) {\bf Fake projective planes arising from $\cT =\{2,3\}$}: Let $\cT =\{2,3\}$. Since  $e'(P_2) = 3$, and $e'(P_3)= 7$, see 3.5 and 2.5({\it ii}),\,({\it iii}), $$\mu(G(\bR)/\Lambda) = \frac{1}{21}e'(P_2)e'(P_3) = 1.$$ Hence, the orbifold Euler-Poincar\'e characteristic $\chi(\Lambda)$ of $\Lambda$ is $3$. As the maximal normal pro-3 subgroup of the non-hyperspecial maximal parahoric subgroup $P_3$ of $G(\bQ_3)$ is of index $96$, $P_3$ does not contain any elements of order $7$.  But any nontrivial element of $G(\bQ)$ of finite order is of order $7$ (Lemma 5.6),  so we conclude that $\Lambda$ ($\subset G(\bQ)\cap P_3$) is torsion-free. As in 5.9, using Theorem 15.3.1 of [Ro] we conclude that $B/\Lambda$ is a fake projective plane, $\Lambda$ ($\cong\overline\Lambda$) is its fundamental group, and since $\overline\Gamma$ is the normalizer of $\Lambda$ in $\overline{G}(\bR)$, $\overline{\Gamma}/\overline{\Lambda}$ is the automorphism group of $B/\Lambda$.  As in 5.10, we 
 see that any torsion-free subgroup $\Pi$ of $\overline\Gamma$ of index $3$ such that $\Pi/[\Pi,\Pi]$ is finite is the fundamental group of a fake projective plane.  

\vskip2mm

\ni (ii) {\bf The fake projective plane arising from $\cT =\{2,5\}$}:  Let $\cT = \{2,5\}$. As $e'(P_5) = 21$, see 2.5({\it iii}), $$\mu(G(\bR)/\Lambda) =\frac{1}{21}e'(P_2)e'(P_5) = 3.$$ Hence, $\chi(\overline\Gamma ) = 3\chi(\Gamma) = 9\mu(G(\bR)/\Gamma) = 9\mu(G(\bR)/\Lambda)/[\Gamma :\Lambda] = 3$. From this we conclude that the only subgroup of $\overline\Gamma$ which can be the fundamental group of a fake projective plane is $\overline\Gamma$ itself. Moreover, as $H^1(\overline\Lambda, \bC)$, and hence also $H^1(\overline\Gamma, \bC)$, are trivial (Theorem 15.3.1 of [Ro]), $B/\overline\Gamma$ is a fake projective plane, and $\overline\Gamma$ is its fundamental group, if and only if, $\overline\Gamma$ is torsion free. 
\vskip1mm

We will now show, using a variant of the argument employed in the proof of Proposition 5.8, that $\overline\Gamma$ is torsion-free. Since the maximal normal pro-5 subgroup of the non-hyperspecial maximal parahoric subgroup $P_5$ of $G(\bQ_5)$ is of index $720$, $P_5$ does not contain any elements of order $7$. This implies that $\Lambda$ ($\subset G(\bQ)\cap P_5$), and hence also $\overline\Lambda$, are  torsion-free since any element of $G(\bQ)$ of finite order is of order $7$ (Lemma 5.6). Now as $\overline\Lambda$ is a normal subgroup of $\overline\Gamma$ of index $3$, we conclude that the order of any nontrivial element of $\overline\Gamma$ of finite order is $3$. 
\vskip1mm

Let $\cG$ be the connected reductive $\bQ$-subgroup of ${\rm GL}_{1,\cD}$, which contains $G$ as a normal subgroup, such that $$\cG(\bQ) = \{ z\in \cD^{\times}\, |\, z\sigma(z)\in \bQ^{\times}\}.$$ Then the center $\cC$ of $\cG$ is $\bQ$-isomorphic to $R_{\ell/\bQ}({\rm GL}_1)$. The adjoint action of $\cG$ on the Lie algebra of $G$ induces a $\bQ$-isomorphism $\cG/\cC\to {\overline G}$. As $H^1(\bQ, \cC) =\{0\}$, the natural homomorphism $\cG(\bQ)\to {\overline G}(\bQ)$ is surjective. Now, if possible, assume that $\overline\Gamma$ contains an element of order $3$. We fix a  $g\in \cG(\bQ)$ whose image $\overline{g}$ in $\overline{G}(\bQ)$ is an element of order $3$ of $\overline\Gamma$. Then $\lambda :={g}^3$ lies in $\ell^{\times}$. Let $a = g\sigma(g)\in \bQ^{\times}$. Then the norm of $\lambda$ (over $\bQ$) is $a^3\in {\bQ^{\times}}^3$. Let $x$ be the unique cube-root of $\lambda$ in the field 
$L=\ell[X]/(X^3-\lambda)$. Then there is an embedding of $L$ in $\cD$ which maps $x$ to $g$. We will view $L$ as a field contained in $\cD$ in terms of this embedding. Then $\sigma(x) = a/x$.  The reduced norm of $x$ is clearly $\lambda$, and the image of $\overline{g}$ in $H^1(\bQ, C)\subset \ell^{\times}/{\ell^{\times}}^3$ is the class of $\lambda$ in $\ell^{\times}/{\ell^{\times}}^3$.    As in the proof of Proposition 5.8 (cf.\:also 5.4), we see that since $g$ stabilizes the collection $(P_q)$, $w(\lambda)\in 3\bZ$ for every normalized valuation $w$ of $\ell =\bQ(\sqrt{-7})$ not lying over $2$.    
\vskip1mm

We assert that the subgroup $\ell^{\bullet}_{\{2\}}$ of $\ell^{\times}$ consisting of elements $z$ whose norm lies in ${\bQ^{\times}}^3$, and $w(z)\in 3\bZ$ for every normalized valuation $w$ of $\ell$ not lying over $2$, equals 
${\ell^{\times}}^3 \cup (1+\sqrt{-7}){\ell^{\times}}^3 \cup (1-\sqrt{-7}){\ell^{\times}}^3$.  Since $\ell$ does not contain a nontrivial cube-root of unity, and its class number is $1$, the subgroup $\ell_3$ of $\ell^{\times}$ consisting of $z\in \ell^{\times}$ such that for {\it every} normalized valuation $w$ of $\ell$, $w(z)\in 3\bZ$ coincides with ${\ell^{\times}}^3$ (see the proof of Proposition 0.12 in [BP]), and $\ell_3$ is of index $3$ in the subgroup $\ell^{\bullet}_{\{2\}}$ (of $\ell^{\times}$).  As $(1+\sqrt{-7})(1-\sqrt{-7}) = 8\, (\in {\ell^{\times}}^3)$, it follows that $1+\sqrt{-7}$ and $1-\sqrt{-7}$ are units at every nonarchimedean place of $\ell$ which does not lie over $2$. Moreover, it is easy to see that  if $v'$ and $v''$ are the two normalized valuations of $\ell$ lying over $2$, then neither $v'(1+\sqrt{-7})$ nor $v''(1+\sqrt{-7})$ is a multiple of $3$. This implies, in particular, that $1\pm \sqrt{-7}\notin {\ell^{\times}}^3$.  From these observations, 
 the above assertion is obvious. Now we note that $1\pm\sqrt{-7}$ is not a cube in $\bQ_5(\sqrt{-7})$ (to see this, it is enough to observe, using a direct computation, that $(1\pm\sqrt{-7})^8\ne 1$ in the residue field of $\bQ_5(\sqrt{-7})$). Since $\lambda\in \ell^{\bullet}_{\{2\}}$ and is not a cube in $\ell^{\times}$, it must lie in the set $(1+\sqrt{-7}){\ell^{\times}}^3\cup (1-\sqrt{-7}){\ell^{\times}}^3$. But no element of this set is a cube in $\bQ_5(\sqrt{-7})$. Hence, $\mathfrak{L} :=L\otimes_{\ell}{\bQ_5(\sqrt{-7})}$ is an unramified field extension of $\bQ_5(\sqrt{-7})$ of degree $3$.  
\vskip1mm

Let $T$ be the centralizer of $g$ in $G$. Then $T$ is a maximal $\bQ$-torus of $G$. Its group of $\bQ$-rational points is $L^{\times}\cap G(\bQ)$. The torus $T$ is anisotropic over $\bQ_5$ and its splitting field over $\bQ_5$ is the unramified cubic extension $\mathfrak{L}$ of $\bQ_5(\sqrt{-7})$. This implies  that any parahoric subgroup of $G(\bQ_5)$ containing $T(\bQ_5)$ is hyperspecial. We conclude from this that $T(\bQ_5)$ is contained in a unique parahoric subgroup of $G(\bQ_5)$, and this parahoric subgroup is hyperspecial. According to the main theorem of [PY], the subset of points fixed by $g$ in the Bruhat-Tits building of $G/\bQ_5$ is the building of $T/\bQ_5$. Since the latter consists of a single point, namely the vertex fixed by the hyperspecial parahoric subgroup containing $T(\bQ_5)$, we infer that $g$ normalizes a unique parahoric subgroup of $G(\bQ_5)$, and this parahoric subgroup is hyperspecial.  As $P_5$ is a non-hyperspecial maximal parahoric subgroup of 
  $G(\bQ_5)$, it cannot be normalized by $g$. Thus we have arrived at a contradiction. This proves that $\overline\Gamma$ is torsion-free. Hence, $B/\overline{\Gamma}$ is a fake projective plane, and $\overline\Gamma$  is its fundamental group. Since the normalizer of $\overline\Gamma$ in $\overline{G}(\bR)$ is $\overline\Gamma$, the automorphism group of $B/\overline{\Gamma}$ is trivial. 
\vskip3mm

\ni{\bf A.4.} It has been pointed out by Cartwright that in Proposition 8.6 of [1] the pairs $\cC_{20}$, $\cC_{26}$ and $\cC_{35}$ have erroneously been omitted, and Steger has pointed out  that the  assertion $\cT = \cT_0$ of Proposition 8.6 is incorrect if $(k,\ell) = \cC_{18}$. The corrected Proposition 8.6 is stated  below in which we have listed the possible $\cT$'s. It turns out  that if $(k,\ell)$ is different from $\cC_{18}$ and 
$\cC_{20}$, then $\cT = \cT_0$. 
\vskip2mm

Proposition 8.6 of [1] should be replaced with the following.
\vskip1mm

\begin{center}
\parbox[t]{11.5cm}{{\bf 8.6.\:Proposition.} {\it Assume that $\cD$ is a cubic division algebra. If the orbifold Euler-Poincar\'e characteristic $\chi(\Gamma)$ of $\Gamma$ is a reciprocal integer, then the pair $(k,\ell)$ must be one of the following nine: $\cC_2$, $\cC_3$, $\cC_{10}$, $\cC_{18}$, $\cC_{20}$, $\cC_{26}$, $\cC_{31}$, $\cC_{35}$ and $\cC_{39}$. Moreover, $\cT_0$ consists of exactly one place $\fv$, and $\cT = \cT_0$ except in the case where $(k,\ell)$ is either $\cC_{18}$ or $\cC_{20}$. Except for the pairs $\cC_3$, $\cC_{18}$ and $\cC_{35}$, $\fv$ is the unique place of $k$ lying over $2$; for $\cC_3$ and $\cC_{35}$, it is the unique place of $k$ lying over $5$, and for $\cC_{18}$ it is the unique place of $k$ lying over $3$.}}
\end{center}
 \vskip1mm

\begin{center}
\parbox[t]{11.5cm} {{\it Description of the possible $\cT$ if the pair $(k,\ell)$ is either $\cC_{18}$ or $\cC_{20}$:    If  $(k,\ell)= \cC_{18} = (\bQ(\sqrt{6}),\bQ(\sqrt{6}, \zeta_3))$, the possibilities are $\cT =\cT_0=\{\fv\}$, and $\cT =\{\fv,\fv_2\}$, where $\fv_2$ is the unique place of $k=\bQ(\sqrt{6})$ lying over $2$.  On the other hand, if $(k,\ell) = \cC_{20} = (\bQ(\sqrt{7}), \bQ(\sqrt{7}, \zeta_4))$, let $\fv'_3$ and $\fv''_3$ be the two places of $k =\bQ(\sqrt{7})$  lying over $3$. Then either $\cT = \cT_0 =\{\fv\}$, or $\cT =\{\fv, \fv'_3\}$, or $\cT = \{\fv, \fv''_3\}$.}}    
\end{center}
\vskip1mm

In the proof of this proposition given in [1], the sentence `` We can show that unless $(k,\ell)$ is one of the six pairs... :" should be replaced with ``We can show that unless (i) $(k,\ell)$ {\it is one of the following nine pairs $\cC_2$, $\cC_3$, $\cC_{10}$, $\cC_{18}$, $\cC_{20}$, $\cC_{26}$, $\cC_{31}$, $\cC_{35}$ and $\cC_{39}$}, (ii) $\cT_0$ and $\cT$ {\it are as in the proposition}, and (iii) $P_v$ {\it is maximal for all $v\in V_f$, except when the pair is $\cC_{18}$}, at least one of the following two assertions will hold:".
\vskip2mm
  
\ni{\bf A.5.} {\bf Twenty eight classes of fake projective planes.} In 5.11 of [1] we have described the fake projective planes arising from the pair $(7,2)$ and $\cT =\cT_0 =\{2\}$. We have shown above in A.3 that the pair $(7,2)$ with $\cT =\{2,3\}$ gives us {\it two} distinct  classes of fake projective planes, and with $\cT=\{2,5\}$ it gives {\it two} more classes (each of these classes consists of just one fake projective plane). We will show below  that the pair $(k,\ell) = \cC_{18}$, with $\cT =\{\fv,\fv_2\}$ gives {\it two} classes of fake projective planes, see A.6 and A.7. In addition, the pair $\cC_{20}$ gives {\it three} classes of fake projective planes (A.9), but $\cC_{26}$ and $\cC_{35}$ give none (A.10). These {\it nine} classes are not among the {\it seventeen} classes of fake projective planes described in Sections 5 and 9 of [1]. Also, as observed in A.2, under the finer classification scheme being used here, the fake projective planes arising from the pairs $(k,\ell) = (\bQ, \bQ(\sqrt{-1}))$ and $(\bQ, \bQ(\sqrt{-2}))$ are placed in six classes rather than four. {\it There are thus a total of  twenty eight classes of fake projective planes.} 
\vskip1mm

\vskip1mm

\ni{A.6--A.10} {\it should be added at the end of 9.3 of [1]}.
\vskip1mm

\ni{\bf A.6.} Let $(k,\ell)$ be any of the following  four pairs: $ \cC_{18}=(\bQ(\sqrt{6}),\bQ(\sqrt{6},\zeta_3))$,  $\cC_{20} = (\bQ(\sqrt{7}), \bQ(\sqrt{7},\zeta_4))$,  $\cC_{26} = (\bQ(\sqrt{15}), \bQ(\sqrt{15}, \zeta_4))$, and $\cC_{35} = (\bQ(\zeta_{20}+\zeta_{20}^{-1}), \bQ(\zeta_{20}))$. 
Note that for each of these four pairs, $D_{\ell} = D_k^2$, so every nonarchimedean place of $k$ is unramified in $\ell$. Let $p=3$ for the first pair, $p=2$ for the next two pairs, and $p =5$ for the last pair, and $\fv$ be the unique place of $k$ lying over $p$. Then $\fv$ splits in $\ell$, and $p$ is the cardinality of the residue field of $k_{\fv}$. Let $v_o$ be a fixed real place of $k$. Now let $\cD$, $\sigma$, $G$, $\overline{G}$ be as in 9.1. We fix a coherent collection $(P_v)_{v\in V_f}$ of parahoric subgroups $P_v$ of $G(k_v)$ such that for  every nonarchimedean place $v$ of $k$, $P_v$ is maximal,  and it is hyperspecial for all $v\ne \fv$. Let $\Lambda = G(k)\cap \prod_{v\in V_f} P_v$, and $\Gamma$ be the normalizer of $\Lambda$ in $G(k_{v_o})$. Let $\overline{\Gamma}$, $\overline{\Lambda}$ be the images of $\Gamma$ and $\Lambda$ in $\overline{G}(k_{v_o})$. 
\vskip1mm

{\it In the rest of this section we will deal exclusively with}  $(k,\ell) =\cC_{18}$. Let $\fv_2$ be the unique place of $k= \bQ(\sqrt{6})$ lying over $2$. Note that $\ell = \bQ(\sqrt{6},\zeta_3)=\bQ(\sqrt{-2},\sqrt{-3})$, the class number of $\ell$ is $1$, and $\ell_{\fv_{2}} := k_{\fv_2}\otimes_k \ell$ is an unramified field extension of $k_{\fv_{2}}$. We fix an Iwahori subgroup $I$ of $G(k_{\fv_2})$, and a {\it non-hyperspecial} {maximal} parahoric subgroup $P$ (of $G(k_{\fv_2})$) containing $I$. Let $\Lambda_P =G(k)\cap P\cap \prod_{v\in V_f-\{ \fv_2\}} P_v$ and $\Lambda_I= \Lambda_P\cap  I$.  Let $\Gamma_I$ and $\Gamma_P$ be the normalizers of $\Lambda_I$ and $\Lambda_P$ respectively in $G(k_{v_o})$. Then ${\Gamma}_I\subset {\Gamma}_P$. Let $\overline{\Lambda}_I$, $\overline{\Gamma}_I$, $\overline{\Lambda}_P$ and $\overline{\Gamma}_P$ be the images of $\Lambda_I$, $\Gamma_I$, $\Lambda_P$ and $\Gamma_P$ respectively  in $\overline{G}(k_{v_o})$. Note that $\overline{\Gamma}_P$ is
  contained in $\overline{G}(k)$, see, for example, [BP, Proposition 1.2].
\vskip1mm

 It follows from the result in 5.4 that $$[\overline{\Gamma}_I:\overline{\Lambda}_I] =[\Gamma_I : \Lambda_I] = 9 =[\Gamma_P:\Lambda_P] = [\overline{\Gamma}_P:\overline{\Lambda}_P].$$ For the pair $(k, \ell) = \cC_{18}$, using the value $\chi(\Lambda) =1$ given in 9.1, and the values $e'(I) = 9$, and $e'(P) = 3$ obtained from 2.5 ({\it iii}), we find that $\chi(\Lambda_I) = 9$ and $\chi(\Lambda_P) = 3$, and hence, $\chi(\overline{\Gamma}_I) = 3$ and $\chi(\overline{\Gamma}_P) = 1$. Furthermore, it follows from Theorem 15.3.1 of [Ro] that $H^1(\overline{\Gamma}_I,\bC)$, and for any subgroup $\Pi$ of $\overline{\Gamma}_P$ containing $\overline{\Lambda}_P$, $H^1(\Pi,\bC)$ vanish. 
We conclude from these observations that $\overline{\Gamma}_I$ is the fundamental group of a fake projective plane if and only if it is torsion-free, and a subgroup $\Pi$ of $\overline{\Gamma}_P$ containing $\overline{\Lambda}_P$ is the fundamental group of a fake projective plane if and only if it is torsion-free and is of index $3$ in $\overline{\Gamma}_P$.  
\vskip2mm

\ni{\bf Proposition.} (i) $\overline{\Gamma}_I$ {\it is torsion-free and hence it is the fundamental group of a fake projective plane.} 
\vskip1mm

(ii) {\it There are three torsion-free subgroups of $\overline{\Gamma}_P$ containing 
$\overline{\Lambda}_P$ which are fundamental groups of fake projective planes.}
    
\vskip1mm

\ni{\it Proof.} Let $\cG$ be the connected reductive $k$-subgroup of ${\rm GL}_{1,\cD}$, which contains $G$ as a normal subgroup, such that $$\cG(k) = \{ z\in \cD^{\times}\, |\, z\sigma(z)\in k^{\times}\}.$$ Then the center $\cC$ of $\cG$ is $k$-isomorphic to $R_{\ell/k}({\rm GL}_1)$. The adjoint action of $\cG$ on the Lie algebra of $G$ induces a $k$-isomorphism $\cG/\cC\to {\overline G}$. As $H^1(k, \cC) =\{0\}$, the natural homomorphism $\cG(k)\to {\overline G}(k)$ is surjective.
\vskip1mm

Let $C$ be the center of $G$, and $\varphi: \,G\rightarrow \overline{G}$ be the natural isogeny. 
Let $\delta:\, \overline{G}(k)\rightarrow H^1(k,C)\subset \ell^{\times}/{\ell^{\times}}^3$ be the 
coboundary homomorphism. Its kernel is $\varphi(G(k))$. 
Given $\overline{g}\in \overline{G}(k)$, let $g$ be any element of $\cG(k)$ which maps onto $\overline{g}$. Then $\delta(\overline{g}) ={\rm Nrd}(g)$ modulo ${\ell^{\times}}^3$. 
\vskip1mm

Since  $\overline{\Lambda}_I$ is torsion-free (cf.\,Lemma 9.2), and $[\overline{\Gamma}_I:\overline{\Lambda}_I] = 9$, if $\overline{\Gamma}_I$ contains an element of finite order, then it contains an element $\overline{g}$ of order $3$. We fix an element $g\in \cG(k)$ which maps onto $\overline{g}$. Then $a :=g\sigma(g)\in k^{\times}$, and $\lambda : = g^3$ lies in $\ell^{\times}$. The reduced norm of $g$ is clearly $\lambda$; the norm of $\lambda$ over $k$ is $a^3 \in {k^{\times}}^3$.   Hence, the image $\delta(\overline{g})$ of $\overline{g}$ in $H^1(k,C)\,(\subset \ell^{\times}/{\ell^{\times}}^3)$ is the class of $\lambda$ in $\ell^{\times}/{\ell^{\times}}^3$.   Since $\overline{g}$ stabilizes the collection $(P_v)_{v\in V_f-\{\fv_2\}}$, as in the proof of Proposition 5.8 (cf.\:also 5.4), we conclude that $w(\lambda)\in 3\bZ$ for any normalized valuation of $\ell$ which does not lie over $2 $ or $3$. But as $\fv_2$ does not split in $\ell$, and the norm of  $\lambda$ lies
  in ${k^{\times}}^3$, it is automatic that for the normalized valuation $w$ of $\ell$ lying over $2$, $w(\lambda)\in 3\bZ$.   Therefore, $\lambda\in \ell^{\bullet}_{\{ 3\}}$, where the latter denotes the subgroup of $\ell^{\times}$ consisting of $z$ such that $N_{\ell/k}(z)\in {k^{\times}}^3$, and for all normalized valuations $w$ of $\ell$, except for the two lying over $3$, $w(z)\in 3\bZ$. Now let  $\alpha= (1+\sqrt{-2})/(1-\sqrt{-2})$. It is not difficult to see that $\ell^{\bullet}_{\{ 3\}} = \bigcup_{0\leqslant m,\, n<3}\  \alpha^m\zeta_3^n {\ell^{\times}}^3$.

\vskip1mm

Let $L$ be the field extension of $\ell$ in $\cD$ generated by $g$. Let $T$ be the centralizer of $g$ in $G$. Then $T$ is a maximal $k$-torus of $G$; its group of $k$-rational points is $L^{\times}\cap G(k)$. It can be shown that if $\lambda =g^3\in \alpha^m\zeta_3^n{\ell^{\times}}^3$, with $0\leqslant m, n <3$, then $\ell_{\fv_{2}}\otimes_{\ell}L$ is the direct sum of three  copies of $\ell_{\fv_{2}}$, each stable under $\sigma$, if $n =0$, and it is an unramifield field extension of $\ell_{\fv_{2}}$ of degree $3$ if $n\ne 0$. We conclude from this that the $k$-torus $T$ is anisotropic over $k_{\fv_{2}}$.

\vskip1mm

According to the main theorem of [PY], the subset of points fixed by $g$ in the Bruhat-Tits building of $G/k_{\fv_{2}}$ is the building of $T/k_{\fv_{2}}$. But as $T$ is anisotropic over $k_{\fv_{2}}$, the building of $T/k_{\fv_{2}}$ consists of a single point. Since the two maximal parahoric subgroups of $G(k_{\fv_{2}})$ containing $I$ are nonisomorphic, if $g$ normalizes $I$, then it fixes the edge corresponding to $I$ in the Bruhat-Tits building of $G/k_{\fv_{2}}$. But as $g$ fixes just a single point in this building, we conclude that $g$ (and hence $\overline{g}$) cannot normalize $I$. This proves that $\overline{\Gamma}_I$ is torsion-free, and we have proved assertion (i) of the proposition.

\vskip1mm

We will now prove assertion (ii) of the proposition. It can be seen, using Proposition 2.9 of [BP], cf.\:5.4, that,  under the homomorphism induced by $\delta$, $\overline{\Gamma}_P/\overline{\Lambda}_P$ is isomorphic to  the subgroup $\ell^{\bullet}_{\{ 3\}}/{\ell^{\times}}^3$ of $\ell^{\times}/{\ell^{\times}}^3$. As has been noted above, $\ell^{\bullet}_{\{ 3\}} = \bigcup_{0\leqslant m, \, n<3}\  \alpha^m\zeta_3^n {\ell^{\times}}^3$, and hence, $\ell^{\bullet}_{\{ 3\}}/{\ell^{\times}}^3$ is isomorphic to $\bZ/3\bZ\times \bZ/3\bZ$. There are three subgroups of $\ell^{\bullet}_{\{ 3\}} /{\ell^{\times}}^3$ of index $3$ generated by an element of the form $\alpha^m\zeta^n_3$ with $n\ne 0$. Let $\Pi$ be the inverse image in $\overline{\Gamma}_P$ of any of these three subgroups. Then, as we will show presently, $\Pi$ is torsion-free and so it is the fundamental group of a fake projective plane.

\vskip1mm

Let us assume that $\Pi$ contains a nontrivial element $\overline{g}$ of finite order. Since $\overline{\Lambda}_P$ is torsion-free (cf.\:Lemma 9.2), and $[\Pi :\overline{\Lambda}_P] = 3$, the order of $\overline{g}$ is $3$.  As in the  proof of assertion (i), we fix $g\in \cG(k)$ which maps onto $\overline{g}$, and let $\lambda = g^3$. Then $\lambda$ is the reduced norm of $g$ and it lies in $\ell^{\bullet}_{\{3\}}$. The image $\delta(\overline{g})$ of $\overline{g}$ in $\ell^{\bullet}_{\{3\}}/{\ell^{\times}}^3$ is the class of $\lambda$ modulo ${\ell^{\times}}^3$. 
Since $\Pi$ is the inverse image in $\overline{\Gamma}_P$ of the subgroup generated by $\alpha^m\zeta^n_3$ for some $m,n<3$, with $n\ne 0$, and $\lambda$ is not a cube in $\ell$, $\lambda\in  (\alpha^m\zeta_3^n){\ell^{\times}}^3\cup  (\alpha^m\zeta_3^n)^2{\ell^{\times}}^3 $. Let $L$ be the field extension of $\ell$ in $\cD$ generated by $g$, and let $T$ be the centralizer of $g$ in $G$. Then $T(k) = L^{\times}\cap G(k)$. As observed in the proof of assertion (i), $T$ is a maximal $k$-torus of $G$ which is anisotropic over $k_{\fv_2}$, and its splitting field over $k_{\fv_2}$ is clearly $\ell_{\fv_2}\otimes_{\ell} L$ which is an unramified field extension of $\ell_{\fv_2}$ of degree $3$.  This implies that the unique point in the Bruhat-Tits building of $G/k_{\fv_2}$ fixed by $g$ is hyperspecial. But since $P$ is a 
non-hyperspecial maximal parahoric subgroup of $G(k_{\fv_2})$, it cannot be normalized by $g$. This implies that $\overline{g}$ does not lie in $\overline{\Gamma}_P$, and we have arrived at a contradiction.

\vskip2mm

\ni{\bf A.7. Remark.}  The above proposition implies that the pair $\cC_{18}$  gives two classes of fake projective planes
with $\cT = \{ \fv, \fv_2\}$:\:the class consisting of a unique fake projective plane with the fundamental group isomorphic to $\overline{\Gamma}_I$, and the class consisting of the fake projective planes whose fundamental group is embeddable in $\overline{\Gamma}_P$, but not in $\overline{\Gamma}_I$. Cartwright and Steger have shown that the latter class consists of just three fake projective planes, the ones with the fundamental group as in (ii) of the above proposition.
 
 \vskip3mm

\ni{\bf A.8.} We shall assume now that $(k,\ell)$ is either $\cC_{20}$, $\cC_{26}$ or $\cC_{35}$.  
Let $p$, the place $\fv$ of $k$, the coherent collection $(P_v)_{v\in V_f}$ of parahoric subgroups $P_v$ of $G(k_v)$, and the subgroups $\Gamma$, $\Lambda$, $\overline{\Gamma}$ and $\overline{\Lambda}$ be as in A.6. We recall (see 9.1) that there exists a cubic division algebra $\fD$ with center $k_{\fv}$ such that $G(k_{\fv})$ is the compact group $\mathrm{SL}_1(\fD)$ of elements of reduced norm $1$ in $\fD$.  The first congruence subgroup $G(k_{\fv})^+ :={\mathrm{SL}}_1^{(1)}(\fD)$ of $G(k_{\fv})=\mathrm{SL}_1(\fD)$ is the unique maximal  normal pro-$p$ subgroup of $G(k_{\fv})$, and the quotient $G(k_{\fv})/G(k_{\fv})^+$ is of order $p^2+p+1$. 
\vskip1mm

In case $(k,\ell) = \cC_{20} = (\bQ(\sqrt{7}), \bQ(\sqrt{7},\zeta_4))$,  let $\fv'_3$ and $\fv''_3$ be the two places of $k= \bQ(\sqrt{7})$ lying over $3$. Note that these places do not split in $\bQ(\sqrt{7}, \zeta_4)$. We fix non-hyperspecial {\it maximal} parahoric subgroups $P'$ and $P''$ of $G(k_{\fv'_3})$ and $G(k_{\fv''_3})$ respectively. In this case (i.e., where $(k,\ell)=\cC_{20}$), let  $\Lambda^+ = \Lambda\cap G(k_{\fv})^+$, $\Lambda' = G(k)\cap P'\cap\prod_{v\in V_f-\{\fv'_3\}}P_v$ and $\Lambda'' = G(k)\cap P''\cap \prod_{v\in V_f-\{\fv''_3\}}P_v$.  Then $\chi(\Lambda') = 3\mu (G(k_{v_o})/\Lambda') = 3 = 3\mu(G(k_{v_o})/\Lambda'')= \chi(\Lambda'')$. By the strong approximation property, $\Lambda^+$ is a subgroup of index $7$ ($= [G(k_{\fv}):G(k_{\fv})^+]$) of $\Lambda$. Hence, $\chi(\Lambda^+) = 3\mu(G(k_{v_o})/\Lambda^+) = 21\mu(G(k_{v_o})/\Lambda) =3$. 

\vskip2mm

\ni{\bf Lemma.} (1) {\it If $(k,\ell) = \cC_{20}$, then any nontrivial element of $G(k)$ of finite order has order $7$; $\Lambda^+$, $\Lambda'$ and $\Lambda''$ are  torsion-free.}
   
\ni (2) {\it  If $(k,\ell)$ is either $\cC_{26}$ or $\cC_{35}$, then $G(k)$ is torsion-free.}
 \vskip2mm
 
\ni {\it Proof.} Let $x\in G(k)$ ($\subset \cD$) be a nontrivial element of finite order, say of order $m$. As the reduced norm of $-1$ is $-1$, $-1\notin G(k)$, and so $m$ is necessarily odd. Let $L$ be the $\ell$-subalgebra of $\cD$ generated by $x$. Then $L$ is a field extension of $\ell$ of degree $3$ since in the cases presently under consideration, $\ell$ does not contain a nontrivial cube root of unity, and hence $x\notin \ell$.  Therefore, $[L:\bQ] = 6[k:\bQ] = 6d$.

We will first prove the first assertion. Let us assume that $(k,\ell) = \cC_{20}=(\bQ(\sqrt{7}),\bQ(\sqrt{7},\zeta_4))$. In this case, $L$ is of degree $12$ over $\bQ$, and as $\zeta_4\in \ell$, $L$ contains a primitive $4m$-th root of unity. This implies that $\phi(4m)$ divides $12$. From this we conclude that $m$ is either $3$, $7$ or $9$. Now since $G(k_{\fv})^+$  is a normal pro-2 subgroup of index $7$ in $G(k_{\fv})$, it is clear that the order of a nontrivial element of $G(k_{\fv})$ of odd order can only be $7$, and moreover, $G(k_{\fv})^+$ does not contain any nontrivial elements of odd order. We observe now that if $P'^+$ and $P''^+$ are the unique maximal normal pro-3 subgroups of $P'$ and $P''$ respectively, then $[P': P'^+] = 2^5\cdot 3=[P'':P''^+]$, and hence any nontrivial element of odd order of either $P'$ or $P''$ is of order $3$. Assertion (1) follows at once from these observations.

We will now prove the second assertion. If $$(k,\ell) = \cC_{26}=(\bQ(\sqrt{15}), \bQ(\sqrt{15},\zeta_4)),$$ then again $L$ is of degree $12$ over $\bQ$, and as $\zeta_4\in\ell$, we conclude, as above, that $m$ is either $3$, $7$ or $9$. As in the case considered above, $G(k_{\fv})^+$  is a normal pro-2 subgroup of index $7$ in $G(k_{\fv})$, therefore the order of any nontrivial element of $G(k_{\fv})$ of odd order can only be $7$. This implies that $\zeta_7\in L$, and hence, $L =\bQ(\zeta_{28})$. Since the only primes which ramify in this field are $2$ and $7$, whereas $3$ ramifies in $k = \bQ(\sqrt{15}) \subset  L$, we conclude that $G(k)$ is torsion-free if $(k,\ell) =\cC_{26}$. 

Let us now consider $(k,\ell) = \cC_{35}=(\bQ(\zeta_{20}+\zeta_{20}^{-1}), \bQ(\zeta_{20}))$. In this case, $L$ is of degree $24$ over $\bQ$, and as  $L$ is an extension of degree $3$ of $\ell$, and $\zeta_{20}$ (and hence $\zeta_5$) lies in the latter, $5$ does not divide $m$. But there does not exist  such an  $m\ne 1$ for which  $\phi(20m)$ divides $24$. This implies that $G(k)$ is torsion-free.   
\vskip3mm

\ni{\bf A.9. Three classes of fake projective planes arising from $(k,\ell)=\cC_{20}$.}  Let $\Lambda$, 
$\Lambda^+$, $\Lambda'$ and $\Lambda''$ be as in A.7. Let $\Gamma$, $\Gamma'$ and $\Gamma''$ be the normalizers of $\Lambda$, $\Lambda'$ and $\Lambda''$ in $G(k_{v_o})$, and $\overline{\Gamma}$, ${\overline{\Gamma}}'$ and ${\overline{\Gamma}}''$ be their images in $\overline{G}(k_{v_o})$. Let 
${\overline{\Lambda}}^+$, ${\overline{\Lambda}}'$ and ${\overline{\Lambda}}''$ be the images of 
$\Lambda^+$, $\Lambda'$ and $\Lambda''$ in $\overline{G}(k_{v_o})$. By the above lemma, these groups are torsion-free.

Theorem 15.3.1 of [Ro] implies that the first cohomology (with coefficients $\bC$) of 
${\overline{\Lambda}}^+$, ${\overline{\Lambda}}'$ and ${\overline{\Lambda}}''$ vanish. As the Euler-Poincar\'e characteristic of each of these three groups is $3$, we conclude that these groups are the fundamental groups of fake projective planes $B/{\overline{\Lambda}}^+$, $B/{\overline{\Lambda}}'$ and $B/{\overline{\Lambda}}''$ respectively. The automorphism groups of these fake projective planes are respectively $\overline{\Gamma}/\overline{\Lambda}^+$, $\overline{\Gamma}'/\overline{\Lambda}' $, and $\overline{\Gamma}''/\overline{\Lambda}''$, which are of order $21$, $3$ and $3$.  Any subgroup $\Pi$ of $\overline{\Gamma}$ (resp., ${\overline{\Gamma}}'$ or ${\overline{\Gamma}}''$) of index $21$ (resp., $3$),  with vanishing $H^1(\Pi, \bC)$, is the fundamental group of a fake projective plane, namely, that of $B/\Pi$. We thus obtain three distinct classes of fake projective planes from $\cC_{20}$.

\vskip2mm

\ni{\bf A.10.} The proof given in 9.7 can be used almost verbatim to show that {\it the pairs $\cC_{26}$ and $\cC_{35}$ cannot give rise to any fake projective planes.} The point to note is that $\chi(\Lambda) =\chi(\overline{\Lambda})=9$, and $[\overline{\Gamma}:\overline{\Lambda}] =3$. So the only subgroup of $\overline\Gamma$ which can possibly be the fundamental group of a fake projective plane is $\overline\Gamma$ itself. But the argument used in 9.7 allows us to show that $\overline\Gamma$ is not torsion-free (it contains an element of order $3$) and hence it cannot be the fundamental group of a fake projective plane.   

Only the following modification in the proof given in 9.7 is required to treat $\cC_{26}$ and $\cC_{35}$. The class number of the field $\ell$ appearing in $\cC_{26}$, namely the field $\bQ(\sqrt{15},\zeta_4)$ is $2$ (the class number of the field $\ell =\bQ(\zeta_{20})$  in $\cC_{35}$ is $1$). Now if the class number of $\ell$ is either 1 or 2, then we can find an element $a\in\ell^{\times}$ such that $v'(a)$ is either 1 or 2, and for all the other normalized valuations $v$ of $\ell$, $v(a) = 0$. We can work without any difficulty with such an element $a$. 

\vskip3mm
\ni{\bf A.11 Remark.} Theorem 10.1 of [1] holds, with the same proof as given in [1], also  for all the fake projective planes belonging to the new classes described above. Therefore, for every fake projective plane $P$, $H_1(P,\bZ)$ is nontrivial. This provides the following  homological characterization of the complex projective plane $\bP_{\bC}^2.$

\vskip3mm

\ni{\bf Theorem.} {\it A smooth complex surface with the same singular homology groups as $\bP_{\bC}^2$ is biholomorphic to $\bP_{\bC}^2.$}

\vskip3mm

\ni{\bf A.12. Remark}. Proposition 10.3 of [1] holds, with the same proof as given in [1], for all the fake projective planes belonging to the new classes described above except for the ones arising from the pair $\cC_{18}$. We conclude as in 10.4 of [1] that the canonical line bundle is three times a holomorphic line bundle for any fake projective plane for which Proposition 10.3 holds. 

\vskip7mm 

\centerline{ {\bf Corrections to} [1]}
\vskip4mm

\ni{\bf A.13.} {\it We will use the finer classification scheme mentioned in A.2. }
\vskip1mm

Sect.\,5.13 and Theorem 5.14 of [1] should be replaced with the following.

\vskip2mm

\parbox[t]{11.5cm}{{\bf 5.13.} We recall that the hyperspecial parahoric subgroups of $G(k_v)$ are conjugate to each other under $\overline{G}(k_v)$, see [Ti2, 2.5]; moreover, if $v$ does not split in $\ell$, and is unramified in $\ell$, then the non-hyperspecial maximal parahoric subgroups of $G(k_v)$ are conjugate to each other under $G(k_v)$. Using the observations in 2.2, and Proposition 5.3, we see that if $(a,p)\ne (15,2)$ (resp.,\,$(a,p)=(15,2)$), then up to conjugation by $\overline{G}(\bQ)$, there are exactly 2 (resp.,\,4) coherent collections  $(P_q)$ of maximal parahoric subgroups such that $P_q$ is hyperspecial whenever $q$ does not ramify in $\bQ(\sqrt{-a})$ and $q\ne p$, since if $a\ne 15$ (resp.,\,$a=15$), there is exactly one prime (resp.,\,there are exactly two primes, namely 3 and 5) which ramify in $\ell = \bQ(\sqrt{-a})$. Moreover, for $(a,p) = (1,5)$ and $(2,3)$, up to conjugation by $\overline{G}(\bQ)$, there is exactly one coherent collection $(P_q)$ of parahoric subgroups such that $P_q$ is hyperspecial for $q\ne 2,\, p$, and $P_2$, $P_p$ are Iwahori subgroups;  for $(a,p) = (7,2)$, if either $\cT=\{2,3\}$ or $\{2,5\}$, then up to conjugation by $\overline{G}(\bQ)$, there are exactly 2 coherent collections $(P_q)$ of maximal parahoric subgroups such 
that $P_q$ is hyperspecial if, and only if, $q\notin \cT\cup \{7\}$.}

\vskip1mm

\begin{center}
\parbox[t]{11.5cm}{\ \ \ From the results in 5.10--5.12, and A.2, A.3 above,  we conclude that for $(a,p) $ equal to either $(1,5)$ or  $(2,3)$, there are two distinct classes with $\cT =\{p\}$, and one more class with $\cT = \{2,p\}$; for $(a,p) =(23, 2)$, there are two distinct classes;  for $(a,p) =(7,2)$, there are six distinct finite classes, and for $(a,p)=(15,2)$, there are four distinct finite classes, of fake projective planes. Thus the following theorem holds.}
\end{center}

\vskip2mm

\begin{center}
\parbox[t]{11.5cm}{{\bf 5.14. Theorem.} {\it There exist exactly eighteen distinct classes of fake projective planes with $k =\bQ$.} }
\end{center}   
\vskip2mm

\ni{\bf A.14.}  The title of Sect.\:9 should be changed to ``Ten additional classes of fake projective planes",  and $$\cC_{20}= \big( \bQ(\sqrt{7}),\bQ(\sqrt{7},\zeta_4)\big),$$ $$\cC_{26} = \big( \bQ(\sqrt{15}),\bQ(\sqrt{15}, \zeta_4)\big),$$ $$\cC_{35} = \big( \bQ(\zeta_{20}+\zeta^{-1}_{20}), \bQ(\zeta_{20})\big)$$ should be added to the list appearing at the beginning of this section and in the table in 9.1. For these three pairs,  $p$ and the place $\fv$ of $k$ should be as in A.6, and the values of $(q_{\fv}, \mu, \chi(\Lambda))$ are $(2, 1/21, 3/7)$, $(2,1,9)$ and $(5, 1/2^5, 9)$ respectively. These values should be included in the table in 9.1. 

In Sect.\:9.6 the word ``five'' appearing in the first sentence should be replaced with ``ten", and in the last sentence, the assertion ``the pair $\cC_{18}$ gives only one" should be replaced by ``the pair $\cC_{18}$ gives three". At the end of 9.6, add the following ``The pair $\cC_{20}$ gives us three additional classes of fake projective planes."

In Theorem 9.8 the word ``{\it five}" occurring in the first and the last sentence should be replaced with ``{\it ten}", and in the last sentence ``{\it the pair $(\bQ(\sqrt{6}),\bQ(\sqrt{6}, \zeta_3))$ gives one more}" should be replaced with ``{\it the pair  $\cC_{18}=(\bQ(\sqrt{6}),\bQ(\sqrt{6}, \zeta_3))$ gives three, and the pair $\cC_{20} =(\bQ(\sqrt{7}),\bQ(\sqrt{7},\zeta_4))$ gives three more}".
\vskip2mm

\ni{\bf A.15.} The proof of Proposition 8.8 in [1] for the case $(k,\ell) =\cC_3 =(\bQ(\sqrt{5}),\bQ(\sqrt{5},\zeta_4))$ requires a correction. Towards the end of this proof, on p.\,359 of [1], we have said that it follows from 8.7 that $\Lambda$, and hence 
$\overline{\Lambda}$ contains an element of order $3$. The proof of the assertion ``{\it Moreover, if for every nonarchimedean place $v$ of $k$ which does not split in $\ell$, $\ell_v := k_v\otimes_k \ell$ contains a primitive cube-root of unity, then $\Lambda^{\mathfrak m}$ contains an element of order $3$}"  in 8.7 is incorrect. But if the pair $(k,\ell)$ is $\cC_3$, $\Lambda$ appearing in the proof of Proposition 8.8 does contain an element of order $3$ as has been shown by Steger. Assuming the existence of such an element, the rest of the proof of Proposition 8.8 is correct. 
\vskip3mm

\ni{\bf A.16.}\,After the above additions and corrections have been made in Theorem 4.4, Sect.\:5,  Proposition 8.6 and  9.3 of [1],  the word  ``seventeen"  occurring in the second  and the third paragraph of Sect.\,1.1, and the first paragraph of 1.5, 
of [1] should be replaced with ``twenty eight"  everywhere. In the fourth paragraph of 1.3, for all $v\in V_f$, including those which ramify in $\ell$, we fix a parahoric subgroup $P_v$ of $G(k_v)$ which is {\it minimal} among the parahoric subgroups of $G(k_v)$ normalized by $\Pi$.  The condition that for $v\in\cR_{\ell}$, $P_v$ is a maximal parahoric subgroup of $G(k_v)$ should be dropped from the first paragraph of 1.4, and the last two paragraphs of 1.4 should be replaced with the following two paragraphs:

\vskip1mm

``It will turn out that for every $v\in V_f$, $P_v$ appearing in the preceding paragraph is a maximal parahoric subgroup of $G(k_v)$ 
except if $(k,\ell )$ is either $(\bQ, \bQ(\sqrt{-1}))$ or $(\bQ, \bQ(\sqrt{-2}))$ or equals $\cC_{18} = (\bQ (\sqrt{6}),\bQ (\sqrt{6},\zeta_3 ))$. 
In particular, if $\Pi$, $\Lambda$, $\Gamma$, and the parahoric subgroups $P_v$ are as in 1.3, then for $v\in V_f$, since $P_v$ was assumed to be minimal among the parahoric subgroups of $G(k_v)$ normalized by $\Pi$, if for a $v$, $P_v$ is maximal, then it is the {\it unique} parahoric subgroup of $G(k_v)$ normalized by $\Pi$."
\vskip1mm

``We will prove that there are exactly {\it twenty eight} distinct $\{k,\ell, G, (P_v)_{{v\in V_f}}\}$, with $\cD\ne \ell$. Each of these twenty eight sets determines a unique principal arithmetic subgroup $\Lambda$     ($=G(k)\cap \prod_{{v\in V_f}} P_v$), which in turn determines a unique arithmetic subgroup $\overline\Gamma$ of $\overline{G}(k_{v_o})$ (recall that $\overline\Gamma$ is the image in $\overline{G}(k_{v_o})$ of the normalizer $\Gamma$ of $\Lambda$ in $G(k_{v_o})$). For eighteen of these twenty eight, $k =\bQ$, see Sect.\,5; and there are two with $k=\bQ(\sqrt{2})$, two with $k=\bQ(\sqrt{5})$, and three each with $k=\bQ(\sqrt{6})$ and $k =\bQ(\sqrt{7})$, see Sect.\,9. The pair $(k,\ell) = (\bQ, \bQ(\sqrt{-1}))$ gives three, the pair $(\bQ,\bQ(\sqrt{-2}))$ gives three, the pair $(\bQ,\bQ(\sqrt{-7}))$ gives six, the pair $(\bQ, \bQ(\sqrt{-15}))$ gives four,  and the pair $(\bQ, \bQ(\sqrt{-23}))$ gives two classes of fake projective planes."
\vskip3mm

\ni{\bf A.17.} Finally, in bound (1) on p.\,329 of [1], ``$L_{\ell \vert k}r(3)$'' should be replaced with ``$L_{\ell \vert k}(3)$''.

\vskip5mm

{\it Acknowledgments.} We are grateful to Donald Cartwright and Tim Steger for pointing out corrections and omissions in Theorem 4.4 and Proposition 8.6 of [1]. We thank Shigeaki Tsuyumine for providing us good bounds for the denominators of $\zeta_k(-1)$ and $L_{\ell |k}(-2)$ which allowed us to recheck their values used in [1]. 

The first author was supported by the NSF (grant DMS-0653512) and the Humboldt Foundation.  The second author was
supported partially by grants from NSA (grant H98230-07-1-0007) and NSF (grant DMS-0758078).
\vskip7mm

\centerline{\bf References}
\vskip3mm

\ni[1] Prasad, G.,\:Yeung, S-K.: Fake projective planes. Invent.\,Math.\,{\bf 168}, 321-370 (2007).
\vskip1.5mm

\ni[BP] Borel, A.,\:Prasad, G.: Finiteness theorems for discrete subgroups of bounded covolume in semisimple groups. Publ.\,Math.\,IHES No.\,{\bf 69}, 119--171 (1989).
\vskip1.5mm

\ni[KK] Kharlamov,\:V., Kulikov, V.: On real structres on rigid surfaces, Izv.\,Math. {\bf 66}, 133-150 (2002).
\vskip1.5mm

\ni[PY] Prasad, G.,\:Yu, J.-K.: On finite group actions on reductive groups and buildings. Invent.\,Math.\,{\bf 147}, 545-560 (2002).
\vskip1.5mm
 
\ni [Ro] Rogawski, J.: Automorphic representations of unitary groups in three variables. Annals of 
Math.\,Studies {\bf 123}, Princeton U.\,Press, Princeton (1990).

\vskip5mm

\ni{\sc University of Michigan, Ann Arbor, MI 48109}

\ni e-mail: gprasad@umich.edu

\vskip3mm
\ni{\sc Purdue University, West Lafayette, IN 47907}

\ni{email: yeung@math.purdue.edu}

\end{document}